\documentclass[12]{article}
\usepackage[a4paper,bindingoffset=0.2in,%
left=1in,right=1in,top=1in,bottom=1in,%
footskip=.25in]{geometry}

\usepackage{datetime}
\settimeformat{ampmtime}
\usdate

\usepackage{authblk}

\newcommand{\nwc}{\newcommand}
\nwc{\draftdate}{\today}
\usepackage[pdftex]{graphicx}
\usepackage{psfrag}
\usepackage{wrapfig}
\usepackage{epsf}
\usepackage{amsmath,amssymb}
\usepackage{url}
\usepackage[mathscr]{euscript}
\usepackage{xcolor}

\newtheorem{thm}{Theorem}[section]
\newtheorem{lemma}[thm]{Lemma}

\newtheorem{claim}[thm]{Claim}
\newtheorem{rema}[thm]{Remark}

\newcommand{\barint}{\hbox{$\int$\kern-0.75\intwidth
		\vrule width 0.5\intwidth height 2.4pt depth -2pt\kern0.25\intwidth}}
\newlength\intwidth
\setbox0=\hbox{$\int$}
\intwidth=\wd0

\newcommand\avint{\hbox{\hbox{$\displaystyle \int$}\hbox{\kern-.9em{$-$}}}}

\newcommand\smavint{\hbox{\hbox{$\int$}\hbox{\kern-.75em{$-$}}}}

\nwc{\st}{^{\mbox{\it st}}}

\nwc{\qref}[1]{(\ref{#1})}
\nwc{\veloc}{v}
\nwc{\rhoc}{\beta}
\nwc{\hl}{\hat{L}}

\def\Xint#1{\mathchoice
	{\XXint\displaystyle\textstyle{#1}}%
	{\XXint\textstyle\scriptstyle{#1}}%
	{\XXint\scriptstyle\scriptscriptstyle{#1}}%
	{\XXint\scriptscriptstyle\scriptscriptstyle{#1}}%
	\!\int}
\def\XXint#1#2#3{{\setbox0=\hbox{$#1{#2#3}{\int}$}
		\vcenter{\hbox{$#2#3$}}\kern-.51\wd0}}

\def\dashint{\Xint-}

\nwc{\intRp}{\int_0^\infty}
\nwc{\aint}{\dashint}
\nwc{\aaint}{\dashint}

\newcommand{\cA}{{\cal A}}

\newcommand{\cE }{{\cal E}}

\newcommand{\cR}{{\cal R}}

\newcommand{\R}{\mathbb R}

\newcommand{\be}{\begin{eqnarray}}
	\newcommand{\ee}{\end{eqnarray}}
\newcommand{\nn}{\nonumber}
\newcommand{\ben}{\begin{eqnarray*}}
	\newcommand{\een}{\end{eqnarray*}}

\allowdisplaybreaks

\title{}
\author{}
\date{}

\numberwithin{equation}{section}
\usepackage{indentfirst}
\usepackage{stmaryrd}
\usepackage{upgreek,textgreek}

\begin{document}
	
\title{Nonnegative Weak Solution to the Degenerate Viscous Cahn-Hilliard Equation}

\author{Toai Luong
		\thanks{corresponding author, email: tluong4@utk.edu }
	}
\affil{Department of Mathematics, University of Tennessee, Knoxville, TN 37996, USA. }

\date{}
\maketitle

\begin{abstract}
	The Cahn--Hilliard equation is a widely used model for describing phase separation processes in a binary mixture.
	In this paper, we investigate the viscous Cahn--Hilliard equation with a degenerate, phase-dependent mobility.
	We define the concept of a weak solution and establish the existence of such a solution by taking limits of solutions to the viscous Cahn--Hilliard equation with positive mobility.
	Additionally, assuming that the initial data is positive, we demonstrate that the weak solution remains nonnegative and is not identically zero. 
	Finally, we prove that the weak solution satisfies an energy dissipation inequality.
\end{abstract}

	
\section{Introduction}
In this paper, we examine the viscous Cahn--Hilliard (VCH) equation, given by
\begin{align}
	\label{VCH1}\partial_t u &= \nabla \cdot (M(u)\nabla\mu), \quad x\in \Omega, t \in (0,\infty) \\
	\label{VCH2}\mu &= -\kappa\Delta u + W'(u) + \alpha \partial_t u.
\end{align}
The Cahn--Hilliard equation was originally introduced to model phase separation in binary alloys \cite{Cahn2,Cahn5}. 
The viscous version of this equation, introduced by Novick--Cohen, serves as an approximation to the original model \cite{NovickCohenVCH}. 
Here, $\Omega$ is a bounded domain in $\R^n \; (n=1,2,3,\dots)$ and  
$u(x,t)$ represents the relative concentration of the two phases as it evolves over time and space. 
The function $W(u)$ is a double-well potential with two equal minima at $u^- < u^+$, corresponding to the two pure phases. 
The parameter $\kappa > 0$ is related to the thickness of the transition region between the phases, 
with $\sqrt{\kappa}$ proportional to this thickness. 
The parameter $\alpha > 0$ represents the mobility tensor, 
and the term $\alpha \partial_t u$ accounts for the viscosity of the system \cite{Gurtin1996}. 
Hence the equation \qref{VCH1}--\qref{VCH2} are known as the viscous Cahn--Hilliard equation. 
When $\alpha = 0$, the standard Cahn--Hilliard equation is recovered.

To ensure the problem is well-posed, we supplement the VCH equation with the initial condition
\begin{align}\label{initial}
	u(x,0) = u_0(x) \quad \text{in } \Omega,
\end{align}
and appropriate boundary conditions. 
When equations \qref{VCH1}--\qref{VCH2} are paired with homogeneous Neumann boundary conditions,homogeneous Neumann boundary conditions 
\begin{align*}
	\partial_{\textbf{\textup{n}}} u = \partial_{\textbf{\textup{n}}} \mu = 0
	\quad \text{on } \partial\Omega,
\end{align*} 
where $\textbf{\textup{n}}$ is the exterior unit normal vector to the boundary $\partial\Omega$, 
or periodic boundary conditions (with $\Omega$ being a periodic cell), 
the model exhibits energy dissipation characterized by the Cahn--Hilliard energy functional
\begin{align}\label{CH-ener}
	\cE (u) &= \int_\Omega \left( \frac{\kappa}{2}|\nabla u|^2 + W(u) \right) dx.
\end{align}
In this paper, we focus on the periodic boundary conditions.

Various forms of the double-well potential $W(u)$ have been employed to model two-phase systems \cite{DaiCH-sln}. 
A thermodynamically relevant potential is the logarithmic function derived from a mean-field model \cite{Cahn2}:
\begin{align}\label{W-log}
	W_{\rm log}(u)=\frac{\theta}{2} \left[(1+u)\ln\left( \frac{1+u}{2} \right)+(1-u)\ln\left( \frac{1-u}{2} \right) \right] + \frac{\theta_c}{2}(1-u^2), \; 0 < \theta < \theta_c.
\end{align}
The logarithmic terms in \qref{W-log} correspond to the entropy of mixing. 
The constants $\theta$ and $\theta_c$ are proportional to the absolute temperature and a critical temperature of the system, respectively.  
The condition $\theta < \theta_c$ ensures that $W_{\rm log}$ has two equal minina at $u^{\pm}=\pm(1-\cR(\theta, \theta_c))$, 
where $\cR(\theta, \theta_c)>0$ and $\lim_{\theta\to 0}\cR(\theta, \theta_c)=0$. 
Although $W_{\rm log}(u)$ is smooth at its minimizers $u^{\pm}$, it is only defined in $(-1,1)$. 

To facilitate numerical simulations and theoretical analysis, 
$W_{\rm log}(u)$ is often approximated by smooth potentials, such as $W(u)=\frac{1}{4}(u^2-1)^2$, 
or, more generally, 
\begin{align}\label{W-smooth}
	W(u)= \gamma(u-u^+)^2(u-u^-)^2, \quad \gamma > 0.
\end{align}
This form of the double-well potential is widely used in phase-field modeling.

Another approach is to use a double-well potential with infinite barriers outside of $[u^-,u^+]$ \cite{Blowey}, for instance,
\begin{align}\label{W-db}
	W_{\rm db}(u)=\begin{cases}
		\frac{1}{2}(u-u^+)(u-u^-)&, \; u^- \leq u \leq u^+, \\
		\infty&, \; \text{otherwise.}
	\end{cases}
\end{align}
This form of $W(u)$ has singularities at $u=u^{\pm}$, presenting additional challenges for numerical and theoretical analysis.

The Cahn–Hilliard equation has been used to model a variety of phenomena, 
including the dynamics of two populations \cite{Cohen}, biomathematical modeling of bacterial films \cite{Klapper}, phase separations in polymers \cite{Saxena}, tumor tissue growth \cite{Cogswell}, and certain thin film problems \cite{Oron, Thiele}. 
It also finds applications in image processing, 
such as image inpainting and segmentation \cite{Bertozzi, Zhu}. 
Moreover, the equation plays a crucial role in phase-field methods \cite{Moelans}, 
where it models the behavior of conserved order variables. 
Significant research has been conducted on the Cahn–Hilliard model, 
as seen in works like \cite{Cahn2, Cahn5, Elliott1987, Elliott-CH-sln, Elliott-VCH1, Elliott-VCH2, DaiCH-asymp1, DaiCH-asymp2, DaiCH-sln, LuongCH-min, LuongCH-sln} and references therein.

\subsection{Main result}
In this paper, we consider the domain $\Omega=(0,2\pi)^n$, 
and impose periodic boundary conditions on its boundary $\partial\Omega$. 
We assume a cutoff mobility function defined as:
\begin{align}\label{M(u)}
	M(u)=\begin{cases}
		u &, \; u > 0,\\
		0 &, \; u \leq 0, 
		\end{cases}
\end{align}
which is degenerate when $u \leq 0$. 
The degeneracy of $M(u)$ introduces significant technical challenges. 
Additionally, we assume the double-well potential 
$W \in C^2(\R, \R)$ and $W$ satisfies the following conditions for all $z\in \R$:
\begin{align}
	\label{grow-1} C_1 |z|^{r+1} - C_2 \le &W(z) \le C_3 |z|^{r+1} + C_4, \\
	\label{grow-2} |&W'(z)| \le C_5 |z|^r + C_6, \\
	\label{grow-3} C_7 |z|^{r-1} - C_8 \le &W''(z) \le C_9 |z|^{r-1} + C_{10},
\end{align}
for some $1\leq r < \infty$ if $n=1,2$ and $1 \leq r \leq \frac{n}{n-2}$ if $n \geq 3$,  
and some constants $C_j>0$ ($j=1,\dots, 10$).

Our analysis follows a framework similar to that in \cite{DaiCH-sln, LuongFCH-sln} and consists of two main steps.  
he first step approximates the degenerate mobility $M(u)$ with a non-degenerate mobility $M_\theta(u)$,  defined for $\theta\in(0,1)$, as:
\begin{align}\label{M-theta}
	M_\theta(u)=\begin{cases}
		u &, \; \text{if } u > \theta,\\
		\theta &, \; \text{if } u \leq \theta. 
		\end{cases}
\end{align}
The positive lower bound of $M_\theta(u)$ allows us to establish the existence of a sufficiently regular weak solution to the VCH equation \qref{VCH1}--\qref{VCH2} with this non-degenerate mobility.

\begin{thm}\label{thm-main1}
	Let $u_0\in H^1(\Omega)$. 
	With the potential $W(u)$ satisfying conditions  \qref{grow-1}--\qref{grow-3} and the mobility  
	$M_\theta(u)$ defined by \qref{M-theta}, for any given constant $T>0$, there exists a function $u_\theta$ that 
	satisfies the following conditions:
	\begin{enumerate}	
		\item[(1)] $u_\theta\in L^\infty(0,T;H^1(\Omega))\cap C([0,T];L^p(\Omega))\cap L^2(0,T;H^3(\Omega))$, 
		where $1\leq p < \infty$ if $n=1,2$ and $1 \leq p < \frac{2n}{n-2}$ if $n \geq 3$.
		
		\item[(2)]  $\partial_tu_\theta\in L^2(0,T;H^1(\Omega))$.
		
		\item[(3)] $u_\theta(x,0)=u_0(x)$ for all $x\in\Omega$.
		
		\item[(4)] $u_\theta$ satisfies the VCH equation \qref{VCH1}--\qref{VCH2} in the following 
		weak sense:
		\begin{align}\label{VCH-w1}
			&\int_0^T \int_\Omega \partial_tu_\theta \phi dxdt \nonumber \\ 
			&= -\int_0^T\int_\Omega M_\theta(u_\theta)\biggl(-\kappa\nabla\Delta u_\theta 
			+ W''(u_\theta)\nabla u_\theta 
			+ \alpha \nabla\partial_tu_\theta \biggr)\cdot\nabla\phi dxdt
		\end{align}
		for all $\phi \in L^2(0,T;H^1(\Omega))$. 				
		
		\item[(5)] For any $t\in [0,T]$, the following energy inequality holds:
		\begin{align} \label{ener-ineq1}
			&\int_\Omega \left(\frac{\kappa}{2}|\nabla u_\theta(x,t)|^2
			+W(u_\theta(x,t))\right)dx  
			+ \alpha \int_0^t\int_\Omega |\partial_tu_\theta(x,\tau)|^2 dxd\tau \nonumber \\ 
			&+ \int_0^t\int_\Omega M_\theta(u_\theta(x,\tau)) 
			|-\kappa \nabla\Delta u_\theta(x,\tau) 
			+ W''(u_\theta(x,\tau))\nabla u_\theta(x,\tau)
			+ \alpha \partial_tu_\theta (x,\tau)|^2dxd\tau 
			\nonumber \\ 
			&\leq\int_\Omega \left(\frac{\kappa}{2}|\nabla u_0|^2+W(u_0)\right)dx.
		\end{align}
				
		\item[(6)] \label{th1-ineq}
		In addition, if $u_0(x)>0$ for all $x\in\Omega$, then for any $0<\theta<1$,
		\begin{align}\label{bnd-neg1}
			\mathop{\textup{ess sup}}_{0\leq t \leq T}\int_\Omega\left|(u_\theta(x,t))_-+\theta\right|^2dx
			\leq C(\theta^2+\theta+\theta^{1/2}),
		\end{align}
		where $(u_\theta)_-=\min\{u_\theta,0\}$, 
		and $C$ is a generic positive constant that may depend on $d,T,\Omega,\eta,u_0$ 
		and $C_j(j=1,2,3,4)$ but not 
		on $\theta$.
	\end{enumerate}
\end{thm}

The second step involves taking the limit of $u_\theta$ as $\theta\to 0$. 
The limiting function $u$ exists and, in the weak sense, solves the VCH equation 
\qref{VCH1}--\qref{VCH2} with the degenerate mobility $M(u)$ defined by \qref{M(u)}.

\begin{thm}\label{thm-main2}
	Let $u_0\in H^1(\Omega)$. 
	With the potential $W$ satisfying \qref{grow-1}--\qref{grow-3} and 
	the mobility  $M(u)$ defined by \qref{M(u)}, 
	for any given constant $T>0$, there exists a function $u$ that satisfies the following conditions:
	
	\begin{enumerate}
		\item[(1)] $u\in L^\infty(0,T;H^1(\Omega))\cap C([0,T];L^p(\Omega))$, 
		where $1\leq p < \infty$ if $n=1,2$ and $1 \leq p < \frac{2n}{n-2}$ if $n \geq 3$.
		
		\item[(2)] $\partial_tu\in L^2(0,T;L^2(\Omega))$.

		\item[(3)] $u(x,0)=u_0$ for all $x\in\Omega$.
		
		\item[(4)] $u$ can be considered as a weak solution to the VCH equation \qref{VCH1}--\qref{VCH2} in the following weak sense:
		\begin{enumerate}
			\item[(a)] Let $P$ be the set where $M(u)$ is not degenerate, defined as:
			\begin{align}
				P:=\{(x,t)\in\Omega_T: u(x,t) > 0 \}.
			\end{align}
			There exists a set $B\subset\Omega_T$ with $|\Omega_T\backslash B|=0$ and a function 
			$\Psi:\Omega_T\to\R^n$ satisfying $\chi_{B\cap P}M(u)\Psi\in 
			L^2(0,T;L^2(\Omega,\R^n))$, 
			where $\chi_{B\cap P}$ is the characteristic function of $B\cap P$, such that
			\begin{align}\label{VCH-w2}
				\int_0^T \int_\Omega \partial_tu \phi dxdt
				=-\int_{B\cap P}M(u)\Psi\cdot\nabla\phi dxdt
			\end{align}
			for all $\phi\in L^2(0,T;H^1(\Omega))$.
			
			\item[(b)] If for some open set $U\subset\Omega$, the generalized derivatives $\nabla\Delta u \in L^q(U_T,\R^n)$ and $\nabla\partial_tu \in L^q(U_T,\R^n)$ for some $q>1$, 
			where $U_T=U\times(0,T)$, then
			\begin{align*}
				\Psi=-\kappa\nabla\Delta u + W''(u)\nabla u  + \alpha \nabla \partial_tu
				\quad in\;U_T.
			\end{align*}
			
			Additionally, for any $t\in [0,T]$, the following energy inequality holds:
			\begin{align}\label{ener-ineq2}
				&\int_\Omega \left(\frac{\kappa}{2}|\nabla u(x,t)|^2+W(u(x,t))\right)dx 
				+ \alpha \int_0^t\int_\Omega |\partial_t u(x,\tau)|^2 dxd\tau \nonumber \\
				&+\int_{\Omega_t\cap B\cap P} M(u(x,\tau))|\Psi(x,\tau)|^2dxd\tau \nonumber \\
				&\leq\int_\Omega \left(\frac{\kappa}{2}|\nabla u_0|^2
				+W(u_0)\right)dx,
			\end{align}		
			where $\Omega_t = \Omega \times (0,t)$.
		\end{enumerate}	
		
		\item[(5)] If $u_0(x)>0$ for all $x\in\Omega$, then $u(x,t)\geq 0$ for all $(x,t)\in\Omega_T$, 
		and $u(x,t)$ is not identically zero in $\Omega_T$.
	\end{enumerate}
\end{thm}

\subsection{Notation}
In this paper, $C$ denotes a generic positive constant that may depend on 
$n,T,\Omega,u_0,\kappa,\alpha,m$ and $C_j(j=1,...,10)$, but is independent of $\theta$. 

The structure of this paper is as follows: 
\begin{itemize}
	\item In Section \ref{pos-mobi}, we prove Theorem~\ref{thm-main1} using the Galerkin approximation method.
	
	\item In Section \ref{dege-mobi}, we prove Theorem~\ref{thm-main2}, which establishes the existence of a weak solution to the VCH equation \qref{VCH1}--\qref{VCH2} with degenerate mobility.
\end{itemize}

\section{Weak solution for the positive mobility case} \label{pos-mobi}
In this section, we prove Theorem~\ref{thm-main1}. 
For simplicity, we consider the domain $\Omega = (0, 2\pi)^n$. 
Fix any $T > 0$, and define $\Omega_T:= \Omega\times(0,T)$. 
Let $\{\phi_j:j=1,2,...\}$ be the normalized eigenfunctions in $L^2(\Omega)$, that is,  $||\phi_j||_{L^2(\Omega)}=1$, 
of the eigenvalue problem
\begin{align*}
	-\Delta u&=\lambda u \quad \mbox{in} \; \Omega
\end{align*}
subject to periodic boundary conditions on $\partial\Omega$. 
The eigenfunctions $\phi_j$ form a complete orthonormal basis for $L^2(\Omega)$ 
and are orthogonal in the $H^k(\Omega)$ scalar product for any $k \geq 1$. 
Without loss of generality, we assume that the first eigenvalue $\lambda_1=0$, 
so $\phi_1\equiv(2\pi)^{-\frac{n}{2}}$.

\subsection{Galerkin approximation} \label{Galerkin}
We consider the Galerkin approximation for the viscous Cahn--Hilliard equation \qref{VCH1}--\qref{VCH2}:
\begin{align}
	u^N(x,t)&=\sum_{j=1}^{N}c^N_j(t)\phi_j(x),\quad \mu^N(x,t)=\sum_{j=1}^{N}d^N_j(t)\phi_j(x),	\\	
	\label{ode1} \int_\Omega\partial_tu^N\phi_jdx&=-\int_\Omega M_\theta(u^N)\nabla\mu^N\cdot\nabla\phi_jdx, \\	
	\label{ode2} \int_\Omega\mu^N\phi_jdx&=\int_\Omega(\kappa \nabla u^N \cdot\nabla\phi_j + W'(u^N) \phi_j + \alpha \partial_t u^N\phi_j )dx, \\	
	\label{ode3} u^N(x,0)&=\sum_{j=1}^{N}\left(\int_\Omega 	u_0\phi_jdx\right)\phi_j(x).
\end{align}
This setup leads to an initial value problem for the system of ordinary differential equations for $c^N_1,...,c^N_N$:
\begin{align}
	\label{ode4}\partial_tc_j^N &= -\sum_{k=1}^{N}d_k^N\int_\Omega M_\theta\left(\sum_{i=1}^{N}c^N_i\phi_i\right) \nabla\phi_k\cdot\nabla\phi_jdx, \\	
	\label{ode5}d_j^N &= \kappa \lambda_j c_j^N + 
	\int_\Omega W'\left(\sum_{k=1}^{N}c^N_k\phi_k\right)\phi_jdx
	+ \alpha\partial_tc_j^N, \\	
	\label{ode6}c_j^N(0) &= \int_\Omega u_0\phi_jdx.
\end{align}
Since the right hand side of \qref{ode4} depends continuously on $c^N_1,...,c^N_N$, 
the initial value problem \qref{ode4}--\qref{ode6} has a local solution. 
Using arguments similar to those in Section 2.1 in \cite{DaiCH-sln}, 
we can derive Lemmas~\ref{ini-lem} and \ref{bnd-lem}, 
which provide a uniform bound for $c^N_1,...,c^N_N$. 
Consequently, the initial value problem \qref{ode4}--\qref{ode6} has a global solution.

\begin{lemma}\label{ini-lem}
	Let $u^N$ be a solution of the system \qref{ode1}--\qref{ode3}, 
	then for any $t\in [0,T]$, we have
	\begin{align}
		\int_\Omega u^N(x,t)dx=\int_\Omega u^N(x,0)dx.
	\end{align}
\end{lemma}

\begin{lemma}\label{bnd-lem}
	For any $N=1,2,3,...$, let $u^N$ be a solution of the system \qref{ode1}--\qref{ode3}, 
	then we have
	\begin{align}
		\label{bnd-uN-H1}||u^N||_{L^\infty(0,T;H^1(\Omega))} &\leq C, \\
		\label{bnd-M-muN} \left\|\sqrt{ M_\theta(u^N)}\nabla\mu^N \right\|_{L^2(\Omega_T)} &\leq C, \\
		\label{bnd-W'-Linf}
		||W'(u^N)||_{L^\infty(0,T;L^2(\Omega))} &\leq C, \\
		\label{bnd-M-theta-Linf}
		||M_\theta(u^N)||_{L^\infty(0,T;L^{\frac{n}{2}}(\Omega))} &\leq C, \\
		\label{bnd-uN_t} ||\partial_tu^N||_{L^2(0,T;L^2(\Omega))} &\leq C.
	\end{align}
\end{lemma}

\subsection{Convergence of $\{u^N\}$ and the existence of a weak solution $u_\theta$} \label{weak-sln-gnl}
Since $H^1(\Omega)\subset\subset L^p(\Omega) \hookrightarrow L^2(\Omega)$, 
where $2 \leq p < \infty$ if $n = 1,2$, 
and $2 \leq p < \frac{2n}{n-2}$ if $n \geq 3$, 
by Aubin--Lions Lemma (see \cite{Simon}), we have
\begin{align*}
	\{f\in 	L^2(0,T;H^1(\Omega)):\partial_tf\in L^2(0,T;L^2(\Omega))\}\subset\subset L^2(0,T;L^p(\Omega)),
\end{align*}
and
\begin{align*}
	\{f\in 	L^\infty(0,T;H^1(\Omega)):\partial_tf\in L^2(0,T;L^2(\Omega))\}\subset\subset C([0,T];L^p(\Omega),
\end{align*}
for the indicated values of $p$. 
Thus, by \qref{bnd-uN-H1} and \qref{bnd-uN_t}, there exist a subsequence of $\{u^N\}$ (not relabeled) and a function $u_\theta\in L^\infty(0,T;H^1(\Omega))\cap C([0,T];L^p(\Omega))$ such that, as $N\to\infty$,
\begin{align}
	\label{uN-w*conv-H1} u^N&\rightharpoonup u_\theta \quad \text{weakly--* in } L^\infty(0,T;H^1(\Omega)), \\
	\label{uN-conv1-CLp} u^N&\to u_\theta \quad \text{strongly in } C([0,T];L^p(\Omega)), \\
	\label{uN-conv1-L2Lp} u^N&\to u_\theta \quad \text{strongly in } L^2(0,T;L^p(\Omega)) \; \text{and a.e. in } \Omega_T, \\
	\label{uN_t-wconv} \partial_tu^N&\rightharpoonup \partial_tu_\theta \quad \text{weakly in } L^2(0,T;L^2(\Omega)),
\end{align}
for the number $p$ indicated. 
We also have the following bound for $u_\theta$ due to \qref{bnd-uN-H1} and \qref{bnd-uN_t}:
\begin{align}
	\label{bnd-u-H1}||u_\theta||_{L^\infty(0,T;H^1(\Omega))} + ||\partial_tu_\theta||_{L^2(0,T;L^2(\Omega))}&\leq C.
\end{align}
Since $W'(u) \leq C(1 + |u|^r)$ for $1 \leq r < \infty$ if $n = 1,2$ and $1 \leq r \leq \frac{n}{n-2}$ if $n \geq 3$, 
by \qref{uN-conv1-L2Lp}, we obtain
\begin{align}\label{W'-uN-conv1-Lq} 
	W'(u^N) \to W'(u_\theta) \quad \text{strongly in } C([0,T];L^q(\Omega))
\end{align}
for $1 \leq q < \infty$ if $n = 1,2$ and $1 \leq q < \frac{2n}{r(n-2)}$ if $n \geq 3$. 
Additionally,by \qref{bnd-W'-Linf}, we obtain that
\begin{align}\label{W'-uN-w*conv1-Linf} 
	W'(u^N) \rightharpoonup W'(u_\theta) \quad \text{weakly--* in } L^\infty(0,T;L^2(\Omega)).
\end{align}
Given that $M_\theta(u^N) \geq \theta$, and using \qref{bnd-M-muN}, we have
\begin{align}\label{bnd-grad-muN}
	||\nabla \mu^N||_{L^2(\Omega_T)} \leq C\theta^{-\frac{1}{2}}.
\end{align}
Taking $j=1$ in \qref{ode2}, and using \qref{bnd-W'-Linf} we obtain
\begin{align}
	\int_\Omega \mu^N dx  =  \int_\Omega (W'(u^N) + \alpha \partial_t u^N) dx.
\end{align}
Applying Poincar\'{e}'s inequality and combining with \qref{bnd-W'-Linf} and \qref{bnd-uN_t}, we get
\begin{align}
	||\mu^N||_{L^2(0,T;H^1(\Omega))} \leq C(\theta^{-\frac{1}{2}} + 1).
\end{align}
Hence, there exist a subsequence of $\{\mu^N\}$ (not relabeled) and a function $\mu_\theta\in L^2(0,T;H^1(\Omega))$ such that
\begin{align}\label{muN-wconv-H1}
	\mu^N \rightharpoonup \mu_\theta \quad \text{weakly in } L^2(0,T;H^1(\Omega)).
\end{align}
Since $M_\theta$ is bounded and continuous, by \qref{uN-conv1-CLp}, \qref{uN-conv1-L2Lp} 
and the Bounded Convergence Theorem, we have
\begin{align}
	\label{sqrtM-uN-conv}	\sqrt{M_\theta(u^N)} &\to \sqrt{M_\theta(u_\theta)} \quad \text{strongly in } C([0,T];L^2(\Omega)).
\end{align}
Combining this with \qref{muN-wconv-H1}, we obtain
\begin{align}
	\sqrt{M_\theta(u^N)}\nabla\mu^N \rightharpoonup \sqrt{M_\theta(u_\theta)}\nabla\mu_\theta \quad \text{weakly in } L^2(0,T;L^1(\Omega,\R^n)).
\end{align}
From \qref{bnd-M-muN}, there exist a subsequence of $\{ \sqrt{M_\theta(u^N)}\nabla\mu^N \}$ (not relabeled) and a function $\beta_\theta \in L^2(\Omega_T;\R^n)$ such that
\begin{align}\label{sqrtM-gradmuN-wconv2}
	\sqrt{M_\theta(u^N)}\nabla\mu^N &\rightharpoonup \beta_\theta \quad \text{weakly in } L^2(0,T;L^2(\Omega,\R^n)).
\end{align}
By the uniqueness of weak limits, we have $\beta_\theta=\sqrt{M_\theta(u_\theta)}\nabla\mu_\theta$. 
We also obtain the following bound from \qref{bnd-M-muN}:
\begin{align}
	\int_{\Omega_T} M_\theta(u_\theta)|\nabla\mu_\theta|^2 dxdt \leq C.
\end{align}
Using a similar argument as in Section 2.3 of \cite{DaiCH-sln}, we derive the weak formulation
\begin{align}\label{VCH-w3}
	\int_0^T\int_\Omega \partial_tu_\theta \xi dxdt 
	= -\int_0^T\int_\Omega M_\theta(u_\theta)\nabla\mu_\theta\cdot\nabla\xi dxdt
\end{align}
for all $\xi\in L^2(0,T;H^1(\Omega))$. 

For any $j=1,2,...$ and any function $a_j(t) \in L^2(0,T)$, from \qref{ode2}, we have
\begin{align}
	\int_0^T\int_\Omega\mu^Na_j(t)\phi_jdxdt = \int_0^T\int_\Omega(\kappa \nabla u^N \cdot a_j(t)\nabla\phi_j + W'(u^N) a_j(t)\phi_j + \alpha \partial_t u^N a_j(t)\phi_j )dxdt.
\end{align}
Since  $a_j(t)\phi(x) \in L^2(0,T;C(\bar{\Omega}))$, 
by \qref{uN-w*conv-H1}, \qref{uN_t-wconv}, \qref{W'-uN-w*conv1-Linf} and \qref{muN-wconv-H1}, 
taking limits as $N \to \infty$, we get
\begin{align}
	\int_0^T\int_\Omega\mu_\theta a_j(t)\phi_jdxdt = \int_0^T\int_\Omega(\kappa \nabla u_\theta \cdot a_j(t)\nabla \phi_j + W'(u_\theta) a_j(t)\phi_j + \alpha \partial_t u_\theta a_j(t)\phi_j )dxdt
\end{align}
for any $j=1,2,...$. 
Then for any function $\phi \in L^2(0,T;H^1(\Omega))$, 
since its Fourier series strongly converges to $\phi$ in $L^2(0,T;H^1(\Omega))$, we have
\begin{align}
	\int_0^T\int_\Omega\mu_\theta\phi dxdt = \int_0^T\int_\Omega(\kappa \nabla u_\theta \cdot \nabla \phi + W'(u_\theta) \phi + \alpha \partial_t u_\theta \phi )dxdt
\end{align}
Since $W'(u) \in L^\infty(0,T;L^2(\Omega))$ and $\mu_\theta \in L^2(0,T;H^1(\Omega))$, 
by regularity theory, $u_\theta \in L^2(0,T;H^2(\Omega))$. 
Then we have
\begin{align}\label{mu-theta}
	\mu_\theta = -\kappa \Delta u_\theta + W'(u_\theta) + \alpha \partial_t u_\theta  \quad \text{a.e. in } \Omega_T.
\end{align}

Following a similar argument as in Section 2.4 of \cite{DaiCH-sln}, we have $\nabla W'(u_\theta) \in L^2(0,T;L^2(\Omega))$, 
which implies that $W'(u_\theta) \in L^2(0,T;H^1(\Omega))$. 
Combining with $\mu_\theta \in L^2(0,T;H^1(\Omega))$, 
by \qref{mu-theta}, we have
\begin{align*}
	-\kappa \Delta u_\theta + \alpha \partial_t u_\theta = \mu_\theta -  W'(u_\theta) \in L^2(0,T;H^1(\Omega)),
\end{align*}
which implies
\begin{align*}
	u_\theta \in L^2(0,T;H^3(\Omega)) \quad \text{and} \quad \partial_t u_\theta \in L^2(0,T;H^1(\Omega))
\end{align*}
by regularity theory. 
Thus, 
\begin{align}\label{grad-mu-theta}
	\nabla \mu_\theta = -\kappa \nabla\Delta u_\theta + W''(u_\theta)\nabla u_\theta + \alpha \nabla\partial_t u_\theta  \quad \text{a.e. in } \Omega_T.
\end{align}
Combining \qref{VCH-w3} with \qref{grad-mu-theta} we get
\begin{align}
	&\int_0^T \int_\Omega \partial_tu_\theta \phi dxdt \nonumber \\
	&= -\int_0^T\int_\Omega M_\theta(u_\theta)\biggl(-\kappa\nabla\Delta u_\theta 
	+ W''(u_\theta)\nabla u_\theta 
	+ \alpha \nabla\partial_tu_\theta \biggr)\cdot\nabla\phi dxdt
\end{align}
for all $\phi \in L^2(0,T;H^1(\Omega))$. 

For the initial data, by \qref{ode3}, we have
\begin{align*}
	u^N(x,0) \to u_0(x) \quad \text{in} \; L^2(\Omega).
\end{align*}
Combining with \qref{uN-conv1-CLp}, we have $u_\theta(x,0)=u_0(x)$ for a.e. $x\in\Omega$.

\subsection{Energy inequality}
Since $u^N$ and $\mu^N$ satisfy the following energy identity
\begin{align}
	&\int_\Omega\left(\frac{\kappa}{2}|\nabla u^N(x,t)|^2+W(u^N(x,t))\right)dx 
	+ \int_{0}^{t}\int_\Omega |\partial_t u^N(x,\tau)|^2dxd\tau
	\\ \nonumber
	& + \int_{0}^{t}\int_\Omega  M_\theta(u^N(x,\tau))|\nabla\mu^N(x,\tau)|^2dxd\tau \nonumber \\
	&=\int_\Omega\left(\frac{\kappa}{2}|\nabla u^N(x,0)|^2+W(u^N(x,0))\right)dx, 
\end{align}
using \qref{uN-w*conv-H1}--\qref{uN_t-wconv}, \qref{sqrtM-gradmuN-wconv2} and \qref{grad-mu-theta}, 
by taking limits as $N \to \infty$, we obtain the energy inequality \qref{ener-ineq1}. 

\subsection{Positive initial data}
In this section, we assume that the initial data $u_0(x)>0$ for all $x\in\Omega$. 
We define the entropy densities $\Phi:(0,\infty)\rightarrow [0,\infty)$ by
\begin{align*}
	\Phi''(u)=\frac{1}{M(u)}, \quad \Phi(1)=\Phi'(1)=0,
\end{align*}
and, for each $0<\theta<1$,
\begin{align*}
	\Phi_\theta''(u)=\frac{1}{M_\theta(u)}, \quad \Phi_\theta(1)=\Phi_\theta'(1)=0, \quad \mbox{and} \quad \Phi_\theta\in C^2(\R),
\end{align*}
Using the definition of $M(u)$ and $M_\theta(u)$ in \qref{M(u)} and \qref{M-theta}, we get
\begin{align}\label{Phi(u)}
	\Phi(u)=u\ln u-u+1 \quad\mbox{for}\; u>0,
\end{align}
and
\begin{align}\label{Phi-theta}
	\Phi_\theta(u)=\left\{\begin{array}{l} u\ln u-u+1, \qquad u>\theta,  \\
		\frac{1}{2\theta}u^2+(\ln\theta-1)u+1-\frac{\theta}{2}, \; u\leq\theta. \end{array} \right.
\end{align}
Moreover,  $\Phi_\theta(u)\geq 0$ for all $u\in\R$, $0\leq\Phi_\theta(u)\leq\Phi(u)$ for all $u>0$ and $\Phi_\theta(u)=\Phi(u)$ for all $u\geq\theta$.

Since $\Phi_\theta''$ is bounded, then $\Phi_\theta'(u_\theta)\in L^2(0,T;H^1(\Omega))$. 
So $\Phi_\theta'(u_\theta)$ is an admissible test function for the equation \qref{VCH-w3}. Hence for any $t\in [0,T]$,
\begin{align}\label{Phi-theta-eqn}
	\int_0^t\left\langle \partial_tu_\theta,\Phi_\theta'(u_\theta)\right\rangle_{(H^2(\Omega))',H^2(\Omega)}d\tau =-\int_0^t\int_\Omega M_\theta(u_\theta)\nabla\mu_\theta\cdot\nabla(\Phi_\theta'(u_\theta)) dxd\tau.
\end{align}

\begin{claim}
For each $\theta \in (0,1)$ and for a.e. $t\in [0,T]$,
\begin{align}\label{Phi-theta-eqn-lhs}
	\int_0^t\left\langle \partial_tu_\theta,\Phi_\theta'(u_\theta)\right\rangle_{(H^2(\Omega))',H^2(\Omega)}d\tau = \int_\Omega\Phi_\theta(u_\theta(t,x))dx - \int_\Omega\Phi_\theta(u_0(x))dx.
\end{align}
\end{claim}

\textit{Proof.}
Fix $\theta \in (0,1)$. 
For any $h>0$, we define
\begin{align}
	u_{\theta,h}(t,x):=\frac{1}{h}\int_{t-h}^ tu_\theta(\tau,x)d\tau.
\end{align}
and we set $u_\theta(t,x):=u_0(x)$ for $t\leq 0$,
Then, applying the Aubins--Lions lemma, 
we obtain that
\begin{align} \label{conv-Phi'-u-theta-h}
	\Phi_\theta'(u_{\theta,h}) \rightarrow \Phi_\theta'(u_\theta) \quad \mathrm{strongly \; in} \; L^2(0,T;H^1(\Omega)) \; \mathrm{as} \; h\rightarrow 0,
\end{align}
up to a subsequence. 
Furthermore, we can show that
\begin{align} \label{conv-dt/du-theta-h}
	\partial_tu_{\theta,h} \rightarrow \partial_tu_\theta \quad \mathrm{strongly \; in} \; L^2(0,T;(H^1(\Omega))') \; \mathrm{as} \; h\rightarrow 0.
\end{align}
Indeed, let $\textup{\textbf{J}}_\theta:=M_\theta(u_\theta)\nabla\mu_\theta$, then for any $\xi\in L^2(0,T;H^1(\Omega))$,
\begin{align}
	&\left| \left\langle \partial_tu_{\theta,h} - \partial_tu_\theta,\xi\right\rangle_{(L^2(0,T;(H^1(\Omega))'),L^2(0,T;H^1(\Omega)))} \right| \nonumber \\
	&= \frac{1}{h}\left| \int_0^T\left\langle \int_{t-h}^t(\partial_tu_\theta(\tau) - \partial_tu_\theta(t))d\tau,\xi\right\rangle_{((H^1(\Omega))',H^1(\Omega))} dt \right| \nonumber \\
	&= 	\frac{1}{h}\left| \int_0^T\left\langle \int_{-h}^0(\partial_tu_\theta(t+s) - \partial_tu_\theta(t))ds,\xi\right\rangle_{((H^1(\Omega))',H^1(\Omega))} dt \right| \nonumber \\
	&\leq \frac{1}{h} \int_{-h}^0 \left| \int_0^T \int_\Omega \nabla\xi\cdot (\textup{\textbf{J}}_\theta(t+s)-\textup{\textbf{J}}_\theta(t))dxdt\right| ds \nonumber \\
	&\leq ||\nabla\xi||_{L^2(0,T;H^1(\Omega))}\sup_{-h\leq s \leq 0}||\textup{\textbf{J}}_\theta(\cdot+s)-\textup{\textbf{J}}_\theta(\cdot)||_{L^2(0,T;H^1(\Omega))}.
\end{align}
Hence,
\begin{align}
	||\partial_tu_{\theta,h} - \partial_tu_\theta||_{L^2(0,T;(H^1(\Omega))')} \leq ||\textup{\textbf{J}}_\theta(\cdot+s)-\textup{\textbf{J}}_\theta(\cdot)||_{L^2(0,T;H^1(\Omega))} \rightarrow 0 \quad \mathrm{as} \; h\rightarrow 0.	
\end{align}
So we obtain \qref{conv-dt/du-theta-h}.

Since $\Phi_\theta'(u_{\theta,h})$ and $\partial_tu_{\theta,h}$ are both in $L^2(\Omega_T)$, we have for a.e. $t\in[0,T]$,
\begin{align}
	\int_0^t\left\langle \partial_tu_{\theta,h},\Phi_\theta'(u_{\theta,h})\right\rangle_{(H^1(\Omega))',H^1(\Omega)}d\tau
	&= \int_0^t \int_\Omega \Phi_\theta'(u_{\theta,h})\partial_tu_{\theta,h} dxd\tau \nonumber \\
	&= \int_\Omega \int_0^t \partial_t \Phi_\theta(u_{\theta,h}(\tau,x)) d\tau dx \nonumber \\
	&=\int_\Omega\Phi_\theta(u_{\theta,h}(t,x))dx - \int_\Omega\Phi_\theta(u_0(x))dx.
\end{align}
Passing to the limit as $h\rightarrow 0$ in this equation, combining with \qref{conv-Phi'-u-theta-h} and \qref{conv-dt/du-theta-h}, we obtain \qref{Phi-theta-eqn-lhs}.

\begin{claim}
Let $(u_\theta)_-:=\min\{u_\theta,0\}$, then for any $0<\theta<1$,
\begin{align}\label{bnd-neg2}
	\mathop{\textup{ess sup}}_{0\leq t \leq T}\int_\Omega |(u_\theta)_-+\theta|^2dx\leq C(\theta^2+\theta+\theta^{1/2}).
\end{align}
\end{claim}

\textit{Proof.}
To prove this, for any $z\leq0$, we rewrite $\Phi_\theta(z)$ as
\begin{align}\label{phi-theta-neg}
	\Phi_\theta(z) =\frac{1}{2\theta}(z+\theta)^2+(\ln\theta-2)z+(1-\theta).
\end{align}
Since $0<\theta<1$, \qref{phi-theta-neg} implies
\begin{align}
	(z+\theta)^2 \leq 2\theta \Phi_\theta(z) \quad \mathrm{for\; all} \;z\leq 0. \nn
\end{align}
Hence for any $t\in[0,T]$, since $\Phi_\theta(z) \geq 0$ for all $z\in\R$, we have
\begin{align}\label{bnd-u-theta-neg}
	\int_\Omega |(u_\theta(x,t))_-+\theta|^2dx 
	\leq& 2 \theta\int_\Omega\Phi_\theta(u_\theta(x,t)_-)dx  \nonumber \\
	\leq & 2\theta \left( \int_{\{u_\theta\leq 0\}}\Phi_\theta(u_\theta(x,t))dx+\int_{\{u_\theta>0\}}\Phi_\theta(0)dx\right) 
	\nonumber \\
	\leq& 2\theta \left( \int_\Omega\Phi_\theta(u_\theta(x,t))dx+\int_\Omega\Phi_\theta(0)dx\right)  \nonumber \\
	\leq& 2\theta \left[ \int_\Omega\Phi_\theta(u_\theta(x,t))dx+\left(1-\frac{\theta}{2} \right)|\Omega|\right]. 
\end{align} 
From \qref{bnd-grad-muN} and \qref{bnd-uN-H1}, we have $\|\nabla\mu_\theta\|_{L^2(\Omega_T)} \leq C/\theta^{1/2}$ 
and $\|\nabla u_\theta\|_{L^2(\Omega_T)} \leq C$. Thus by \qref{Phi-theta-eqn}, \qref{Phi-theta-eqn-lhs} and H\"{o}lder's inequality, we have
\begin{align}\
	&\left|\int_\Omega\Phi_\theta(u_\theta(x,t))dx\right| \nn\\
	\leq& \left|\int_\Omega\Phi_\theta(u_0)dx\right| 
	+ \left|\int_0^t\int_\Omega \nabla\mu_\theta(x,\tau)\cdot\nabla u_\theta(x,\tau) dxd\tau\right|  \nonumber\\
	\leq & \int_\Omega\Phi_\theta(u_0)dx + \|\nabla\mu_\theta\|_{L^2(\Omega_T)}\|\nabla u_\theta\|_{L^2(\Omega_T)} 	
	\nonumber \\
	\leq & \int_\Omega\Phi(u_0)dx + \frac{C}{\theta^{1/2}} \nn
\end{align}
for any $t\in[0,T]$. 
This inequality and \qref{bnd-u-theta-neg} implies \qref{bnd-neg2}.

\section{Weak solutions for the degenerate mobility case}\label{dege-mobi}
In this section, we prove Theorem~\ref{thm-main2}, which is the main result of this paper. 
We consider the VCH equation \qref{VCH1}--\qref{VCH2} with the degenerate mobility $M(u)$ defined by \qref{M(u)}.

Fix $u_0\in H^1(\Omega)$ and a sequence of positive numbers $\{\theta_i\}_{i=1}^\infty$ 
that monotonically decreases to 0 as $i\to\infty$. 
For each $\theta_i$, by Theorem~\ref{thm-main1}, there exists a function
\begin{align*}
	u_{\theta_i}\in L^\infty(0,T;H^1(\Omega))\cap C([0,T];L^p(\Omega))\cap L^2(0,T;H^3(\Omega)), 
\end{align*}
whose weak derivative is
\begin{align*}
	\partial_tu_{\theta_i}\in L^2(0,T;L^2(\Omega)),
\end{align*}
where $1\leq p < \infty$ if $n=1,2$ and $1 \leq p < \frac{2n}{n-2}$ if $n \geq 3$, 
such that
\begin{align}
	\label{seq-VCH-w} \int_0^T\int_\Omega \partial_tu_{\theta_i}\xi dxdt &=-\int_0^T\int_\Omega M_{\theta_i}(u_{\theta_i})\nabla\mu_{\theta_i}\cdot\nabla\xi dxdt \quad \text{for any } \xi\in L^2(0,T;H^1(\Omega)),   \\
	\label{seq-mui} \mu_{\theta_i}&=-\Delta u_{\theta_i}+W'(u_{\theta_i})+\alpha\partial_t u_{\theta_i} \quad \text{a.e. in } \Omega_T, \\
	\label{ui(0)} u_{\theta_i}(x,0)&=u_0(x) \quad \text{for all } x\in\Omega.
\end{align}
For simplicity, we denote $u_i:=u_{\theta_i},\mu_i:=\mu_{\theta_i},\omega_i:=\omega_{\theta_i}$, and $M_i:=M_{\theta_i}$.

\subsection{Convergence of $\{u_i\}$ and the existence of a weak solution $u$}
By the proof of Lemma~\ref{bnd-lem}, the bounds on the right hand side of \qref{bnd-uN-H1}--\qref{bnd-uN_t} 
may depend on $n,T,u_0,\Omega,\kappa,\alpha,m$ and $C_j(j=1,...,10)$ but not on $\theta$. 
Thus, there exists a constant $C>0$ independent of $\{\theta_i\}_{i=1}^\infty$ such that
\begin{align}
	\label{bnd-ui-H1}||u_i||_{L^\infty(0,T;H^1(\Omega))} &\leq C, \\
	\label{bnd-M-ui-grad-mui}||\sqrt{M_i(u_i)}\nabla\mu_i||_{L^2(\Omega_T))} &\leq C, \\
	\label{bnd-W'-ui}
	||W'(u^N)||_{L^\infty(0,T;L^2(\Omega))} &\leq C, \\
	\label{bnd-M-ui}
	||M_\theta(u^N)||_{L^\infty(0,T;L^{\frac{n}{2}}(\Omega))} &\leq C, \\
	\label{bnd-ui_t}||\partial_tu_i||_{L^2(0,T;L^2(\Omega))} &\leq C.
\end{align}
Similar to the proof of Theorem~\ref{thm-main1}, 
the boundedness of $\{ u_i \}$ and $\{ \partial_tu_i \}$, 
allows us to extract a subsequence of $\{u_i\}$ (not relabeled) and a function 
\begin{align*}
	u\in L^\infty(0,T;H^1(\Omega))\cap C([0,T];L^p(\Omega))
\end{align*}
with weak derivative
\begin{align*}
	\partial_t u\in L^2(0,T;L^2(\Omega))
\end{align*}
such that, as $i\to\infty$,
\begin{align}
	\label{ui-w*conv-H1}u_i&\rightharpoonup u \quad \text{weakly--* in } L^\infty(0,T;H^1(\Omega)), \\
	\label{ui-conv-C}u_i&\to u \quad \text{strongly in } C([0,T];L^p(\Omega)), \\
	\label{ui-conv-W2p}u_i&\to u \quad \text{strongly in } L^2(0,T;L^p(\Omega)) \; \text{and a.e. in } \; \Omega_T, \\
	\label{ui_t-wconv}\partial_tu_i&\rightharpoonup \partial_tu \quad \text{weakly in } L^2(0,T;L^2(\Omega)),
\end{align}
for any $1\leq p < \infty$ if $n=1,2$ and $1 \leq p < \frac{2n}{n-2}$ if $n \geq 3$. 
By \qref{ui-conv-C} and \qref{ui-conv-W2p} and the uniform convergence of $M_i \to M$ as $i \to \infty$, we have
\begin{align}
	M_i(u_i) \to M(u) \quad \text{a.e. in } \; \Omega_T.
\end{align}

By \qref{bnd-M-ui-grad-mui}, there exists a function $\beta\in L^2(\Omega_T;\R^n)$ such that
\begin{align}
	\label{sqrtMi-grad-mui-conv}\sqrt{M_i(u_i)}\nabla\mu_i\rightharpoonup \beta \quad \text{weakly in } L^2(\Omega_T,\R^n).
\end{align}
Since $M_i(u_i) \to M(u)$ a.e. in $\Omega_T$ and $M_i$ is uniform bounded, 
using a similar argument as in Section 2.3 in \cite{DaiCH-sln}, we get
\begin{align}
	\label{Mi-grad-mui-conv}M_i(u_i)\nabla\mu_i\rightharpoonup \sqrt{M(u)}\beta \quad \text{weakly in } L^2(0,T;L^2(\Omega,\R^n)).
\end{align}
So taking the limit as $i\to\infty$ in \qref{seq-VCH-w}, we obtain
\begin{align}\label{VCH-w3a}
	\int_0^T\int_\Omega \partial_tu\xi dxdt=-\int_0^T\int_\Omega \sqrt{M(u)}\beta\cdot\nabla\xi dxdt
\end{align}
for all $\xi\in L^2(0,T;H^1(\Omega))$.

Now, choose a sequence of positive numbers $\{\delta_j\}_{j=0}^\infty$ that monotonically decreases to $0$. 
For each $\delta_j$, by \qref{ui-conv-W2p} and Egorov's theorem, there exists a subset $B_j\subset\Omega_T$ with $|\Omega_T\backslash B_j|<\delta_j$ such that
\begin{align}\label{ui-conv-unif}
	u_i\to u \quad \text{uniformly in } B_j \; \text{as } i\to\infty.
\end{align}
We can also assume that
\begin{align}\label{B-seq}
	B_1\subset B_2\subset ...\subset B_j\subset B_{j+1}\subset ...\subset\Omega_T.
\end{align}
Define $B:=\cup_{j=1}^\infty B_j$, then $|\Omega_T\backslash B|=0$. 
For each $j$, define the set
\begin{align*}
	P_j:=\{(x,t)\in\Omega_T: u(x,t) > \delta_j\},
\end{align*}
then $P:=\cup_{j=1}^{\infty} P_j = \{(x,t)\in\Omega_T: u(x,t) > 0\}$
is the set where $M(u)$ is not degenerate. 
Using a diagonal argument similar to Section 3.3.2 in \cite{DaiCH-sln}, we are able to find a subsequence 
$\{u_{k,N_k}\}_{k=1}^\infty$, which converges to $u$, 
and a function $\Psi :\Omega_T \to \R^n$ satisfying $\chi_{B\cap P}M(u)\Psi\in L^2(0,T;L^2(\Omega))$ 
such that $u$ and $\Psi$ satisfy the following weak formulation
\begin{align}\label{VCH-w4}
	\int_0^T\int_\Omega \partial_tu \xi dxdt=-\int_{B\cap P}M(u)\Psi\cdot\nabla\xi dxdt
\end{align}
for all $\xi\in L^2(0,T;H^1(\Omega))$. 
Moreover, we can also find a corresponding subsequence $\{\nabla\mu_{k,N_k}\}_{k=1}^\infty$ such that for each $j=1,2,...$,
\begin{align}
	\label{grad-muk-wconv1}	\nabla\mu_{k,N_k} &\rightharpoonup\Psi \quad \text{weakly in } L^2(B_j\cap P_j,\R^n),\\
	\label{Mk-grad-muk-wconv}	\chi_{B_j\cap P_j}\sqrt{M_{k,N_k}(u_{k,N_k})}\nabla\mu_{k,N_k} &\rightharpoonup \chi_{B_j\cap P_j}\sqrt{M(u)}\Psi \quad \text{weakly in } L^2(0,T;L^2(\Omega,\R^n))
\end{align}
as $k\to\infty$. 
As for the initial data, by \qref{ui(0)} and \qref{ui-conv-C}, we have $u(x,0)=u_0(x)$ for all $x\in\Omega$.

\subsection{The relation between $\Psi$ and $u$}
Now we analyze the relation between $\Psi$ and $u$. 
Since
\begin{align*}
	\nabla\mu_i = -\kappa\nabla\Delta u_i + W''(u_i)\nabla u_i + \alpha\nabla\partial_t u_i,
\end{align*}
the desired relation between $\Psi$ and $u$ is
\begin{align*}
	\Psi = -\kappa\nabla\Delta u + W''(u)\nabla u + \alpha\nabla\partial_t u.
\end{align*}
However, due to the regularity of $u$, the terms $\nabla\Delta u$ and $\nabla\partial_t u$ are only defined in the sense of distributions 
and may not be functions.
herefore, we need some higher regularity conditions on $u$.

\begin{claim}\label{claim3}
If $(B_j \cap P_j)^\circ \neq \emptyset$ for some j, where $(B_j \cap P_j)^\circ$ is the interior of $B_j \cap P_j$, then
\begin{align*}
	u \in W^{3,\frac{2n}{n+2}}((B_j \cap P_j)^\circ) \quad \text{and} \quad \partial_tu \in W^{1,\frac{2n}{n+2}}((B_j \cap P_j)^\circ),
\end{align*}
and
\begin{align*}
	\Psi = -\kappa\nabla\Delta u + W''(u)\nabla u + \alpha\nabla\partial_t u \quad \text{in }  (B_j \cap P_j)^\circ.
\end{align*}
\end{claim}

\textit{Proof.}
Using a similar argument as in Section 3.2.1 of \cite{DaiCH-sln}, we obtain
\begin{align}\label{W''-grad-u-conv}
	W''(u_i)\nabla u_i \rightharpoonup W''(u)\nabla u \quad \text{weakly in } L^2(0,T;L^\frac{2n}{n+2}(\Omega,\R^n)) \text{ as } i \to \infty.
\end{align}
 By \qref{grad-mu-theta}, we have, in $(B_j \cap P_j)^\circ$,
\begin{align}
	\nabla(-\kappa\Delta u_{k,N_k} + \alpha \partial_t u_{k,N_k}) = \nabla \mu_{k,N_k} - W''(u_{k,N_k})\nabla u_{k,N_k}.
\end{align}
By \qref{grad-muk-wconv1} and \qref{W''-grad-u-conv}, taking the limits as $k \to \infty$ in the sense of distribution, we get
\begin{align}
	\nabla(-\kappa\Delta u + \alpha \partial_t u) = \Psi - W''(u)\nabla u \quad \text{in } (B_j \cap P_j)^\circ.
\end{align}
Since $\Psi - W''(u)\nabla u \in L^{\frac{2n}{n+2}}(B_j \cap P_j,\R^n) $, 
it follows that
\begin{align*}
	u \in W^{3,\frac{2n}{n+2}}((B_j \cap P_j)^\circ) \quad \text{and} \quad \partial_tu \in W^{1,\frac{2n}{n+2}}((B_j \cap P_j)^\circ),
\end{align*}
and consequently,
\begin{align*}
	\Psi = -\kappa\nabla\Delta u + W''(u)\nabla u + \alpha\nabla\partial_t u \quad \text{in } (B_j \cap P_j)^\circ.
\end{align*}
Claim~\ref{claim3} has been established.

\begin{rema}
The function $\Psi$ is not defined in the set $\Omega_T\backslash(B\cap P)$. 
However, the value of $\Psi$ outside of $B\cap P$ does not affect the integral on the right hand side of \qref{VCH-w4}. 
This ambiguity can be resolved by defining $\Psi$ in any subset of $\Omega_T$ 
where $u$ and $\partial_t u$ have sufficient regularity.
\end{rema}

\begin{claim}\label{claim4}
For any open set $U\subset\Omega$ such that $\nabla\Delta u\in L^q(U_T)$ and $\partial_t u\in L^q(U_T)$ for some $q>1$ that may depend on $U$, where $U_T=U\times(0,T)$, we have
\begin{align}
	\Psi = -\kappa\nabla\Delta u + W''(u)\nabla u + \alpha\nabla\partial_t u \quad \text{in } U_T.
\end{align}
\end{claim}

\textit{Proof.}
Cnsider the limit as $k \to \infty$ of 
\begin{align}\label{grad-muk}
	\nabla\mu_{k,N_k}=-\kappa\nabla\Delta u_{k,N_k} + W''(u_{k,N_k})\nabla u_{k,N_k} + \alpha\nabla\partial_t u_{k,N_k}.
\end{align}
The right hand side weakly converges to $-\kappa\nabla\Delta u + W''(u)\nabla u + \alpha\nabla\partial_t u$ in $L^s(U_T)$, 
where $s = \min \{q , \frac{2n}{n+2}\} >1$. 
Hence,
\begin{align}
	\nabla\mu_{k,N_k}\rightharpoonup -\kappa\nabla\Delta u + W''(u)\nabla u + \alpha\nabla\partial_t u \quad \text{weakly in } L^s(U_T).
\end{align}
Combining this with \qref{grad-muk-wconv1}, we get
\begin{align*}
	\Psi = -\kappa\nabla\Delta u + W''(u)\nabla u + \alpha\nabla\partial_t u \quad\text{in } B\cap P\cap U_T.
\end{align*}
We may extend the function $\Psi$ from $B\cap P\cap U_T$ to $U_T$ 
by defining $\Psi := -\kappa\nabla\Delta u + W''(u)\nabla u + \alpha\nabla\partial_t u$ in $U_T\backslash(B\cap P)$. 
This completes the proof of Claim~\ref{claim4}. 

\

Define the set
\begin{align*}
	\cA:=\bigcup\{U_T=U\times(0,T):\; &U \; \text{is open in }\Omega,\; \nabla\Delta u\in L^q(U_T,\R^n) \text{ and } \nabla\partial_t u\in L^q(U_T,\R^n) \; \\
	& \text{for some } q>1 \;\text{that may depend on } U\},
\end{align*}
then $\cA$ is open in $\Omega_T$ and
\begin{align*}
	\Psi =  -\kappa\nabla\Delta u + W''(u)\nabla u + \alpha\nabla\partial_t u \quad\text{in } \cA.
\end{align*}
So. $\Psi$ is defined in $(B\cap P)\cup\cA$. To extend $\Psi$ to $\Omega_T$, note that
\begin{align*}
	\Omega_T\backslash((B\cap P)\cup\cA)\subset \Omega_T\backslash(B\cap P)=(\Omega_T\backslash B)\cup(\Omega_T\backslash P).
\end{align*}
Since $|\Omega_T\backslash B|=0$ and $M(u)=0$ in $\Omega_T\backslash P$, 
the value of $\Psi$ outside of $(B\cap P)\cup\cA$ does not affect the the integral on the right hand side of \qref{VCH-w4}. 
Thus, we may define $\Psi:=0$ in the complement of $(B\cap P)\cup\cA$.

\subsection{Energy inequality}
By the energy inequality \qref{ener-ineq1}, for any $t \in [0,T]$ and any $j=1,2,...$, we have
\begin{align}\label{ener-ineq3}
	&\int_\Omega \left(\frac{\kappa}{2}|\nabla u_{k,N_k}(x,t)|^2+W(u_{k,N_k}(x,t))\right)dx  
	+ \alpha \int_0^t\int_\Omega |\partial_t u_{k,N_k}(x,\tau)|^2 dxd\tau \nonumber \\
	&+\int_{\Omega_t\cap B_j\cap P_j} M_{k,N_k}(u_{k,N_k}(x,\tau))|\nabla\mu_{k,N_k}(x,\tau)|^2dxd\tau \nonumber \\
	&\leq\int_\Omega \left(\frac{\kappa}{2}|\nabla u_0|^2+W(u_0)\right)dx
\end{align}
Using \qref{ui-w*conv-H1}--\qref{ui_t-wconv} and \qref{Mk-grad-muk-wconv}, 
taking limits as $k\to\infty$ first and then $j\to\infty$ in \qref{ener-ineq3}, we obtain the energy inequality \qref{ener-ineq2}. 

\subsection{Nonnegative weak solution with positive initial data}
Assume the initial data satisfies $u_0(x)>0$ for all $x\in\Omega$. 
From \qref{bnd-neg1}, there exists a constant $C$, independent of $\{\theta_i\}_{i=1}^\infty \subset (0,1)$, such that for each $i=1,2,...$,
\begin{align}\label{bnd-neg3}
	\mathop{\textup{ess sup}}_{0\leq t \leq T}\int_\Omega|(u_i(x,t))_-+\theta_i|^2dx\leq C(\theta_i^2+\theta_i+\theta_i^{1/2}).
\end{align}
Using the convergence results from \qref{ui-conv-C} and \qref{ui-conv-W2p}, and noting that $u$ is continuous in $\Omega_T$, 
we pass to the limit as $i\rightarrow\infty$ in the above inequality.
This yields 
\begin{align}
	\mathop{\textup{ess sup}}_{0\leq t \leq T}\int_\Omega|(u_i(x,t))_-|^2dx\leq 0.
\end{align}
The inequality implies $(u_i(x,t))_- = 0$ a.e. in $\Omega_T$. 
Thus, $u(x,t)\geq 0$ for all $(x,t)\in\Omega_T$. 
Additionally, given that $u_0>0$ in $\Omega$ and considering the continuity of $u$ in $\Omega_T$, 
$u$ is not identically zero in $\Omega_T$. 
This completes the proof of Theorem \ref{thm-main2}.

\section{Conclusion}
In this paper, we have established the existence of a nonnegative weak solution to the degenerate viscous Cahn--Hilliard equation \qref{VCH1}--\qref{VCH2}. 
Understanding the qualitative behavior of these weak solutions is equally important. 
For example, it would be interesting to explore the asymptotic behavior of solutions as certain parameters approach zero. 
In \cite{Kagawa2022}, Kagawa and Ōtani demonstrated that solutions of the VCH equation with constant mobility and homogeneous Dirichlet boundary conditions converge to solutions of the Allen--Cahn equation, the Cahn--Hilliard equation, and the viscous diffusion equation when appropriate parameters are taken to zero. 
Investigating the asymptotic limits of the VCH equation with degenerate mobility presents a challenging yet fascinating problem.

\section{Conflict of Interest}
The author declares that he has no known conflicts of interest that are relevant to the content of this article.

\bibliographystyle{plain}
\bibliography{Luong-bib}

\end{document}